\documentclass[12pt,a4paper]{article}
\usepackage{amsmath}
\usepackage{amssymb}
\usepackage{amsthm}
\usepackage{latexsym}
\usepackage[english,russian]{babel}

\newcommand{\End}{\mathop{\rm End}\nolimits}

\numberwithin{equation}{section} \swapnumbers

\theoremstyle{definition}

\theoremstyle{remark}

\def\?{\marginpar{$\bullet\bullet\bullet$}}

\title{О построении векторных полей сфер}
\author{Игорь~Баяк}
\date{20 февраля 2007}

\begin{document}

\maketitle
\begin{abstract}

В данной работе показано как с помощью группы зеркальных симметрий
сформировать максимально возможную систему линейно независимых
тривиальных векторных полей нечетномерной сферы произвольной
размерности.

\end{abstract}

\section{Постановка задачи и определения}

Необходимая нам группа зеркальных симметрий возникает при решении
конкретной задачи, поэтому сформулируем сначала проблему. Хорошо
известно, что тривиальных векторных полей четномерных сфер не
существует, а максимальное число регулярных линейно независимых
векторных полей нечетномерной сферы $S^{n-1}$ равно
$p(m)=8d+2^{c}-1$, где $c=|m|\mod 4$, $m=4d+c$, $n=(2a+1)2^{m}$,
$m,a\in\mathbb{N}$. Окончательное решение этой задачи нашло
отражение в работе Адамса \cite{ADA}. В свою очередь, Аминов в
монографии \cite{AMI} показал как сформировать, т.е. получить в
явном виде, полную систему регулярных линейно независимых
векторных полей нечетномерных сфер малой размерности. Попробуем
формализовать и обобщить эмпирический результат, полученный
Аминовым.

Пусть в евклидовом пространстве $\mathbb{R}^{n}$ размерности
$n=2^{m}$ уравнением $x^{2}=1$ задана сфера $S^{n-1}$, на которой
лежит произвольный вектор $X$. Ортогональное к $X$ векторное поле
$AX$, где $A$ -- опреатор пространства $\End \mathbb{R}^{n}$,
называется регулярным векторным полем сферы $S^{n-1}$. Поскольку
$(AX,X)=0$, то $A$ это кососимметричный оператор, и поэтому в
матричном представлении $A^{T}=-A$. Перед нами стоит задача ---
найти и представить линейно независимую систему регулярных
векторных полей. Операторы, формирующие эти поля, мы будем искать
среди ортогональных операторов, а условие линейной независимости
заменим условием попарной ортогональности. Тогда, если $A\in
O(n)$, то кососимметричность ортогонального оператора $A$
эквивалентна условию $A^{-1}=-A$ или $A^{2}=-E$. В то же время,
если $A,A'\in O(n)$, $A^{2}=(A')^{2}=-E$ и выполняется условие
ортогональности полей $AX$ и $A'X$, то в силу тождества
$(A'X,AX)=(A^{*}A'X,X)=-(AA'X,X)=0$ мы получим условие
$(AA')^{2}=-E$, из которого следует равносильное ему условие
антикоммутативности операторов $A'A=-AA'$. Итак, нам необходимо
найти максимально возможное число ортогональных кососимметричных
попарно антикоммутирующих операторов, или иначе, ортогональных
кососимметричных операторов, всякое попарное произведение которых
также кососимметрично.

Прежде чем взяться за решение поставленой задачи, напомним вам,
что доказательство существования тривиальных векторных полей
нечетномерных сфер основано на существовании гомотопного
тождественному антиподального отображения нечетномерных сфер.
Обратно, отсутствие векторных полей четномерных сфер связано с
отсутствием гомотопного тождественному антиподального отображения
четномерных сфер. Вместе с тем, всякая сфера $S^{n-1}$ гомеоморфна
границе $n$-мерного параллелепипеда, и поэтому существование
(отсутствие) тривиальных векторных полей сферы $S^{n-1}$
эквивалентно существованию (отсутствию) гомотопного тождественному
антиподального отображения границы $n$-мерного параллелепипеда.
Заметим при этом, что группа изоморфизмов границы $n$-мерного
параллелепипеда равна группе мономиальных подстановок базиса
пространства $\mathbb{R}^{n}$, т.е. она состоит из матриц, в
каждой строке и каждом столбце которых по одному ненулевому
элементу равному $1$ или $-1$, и поэтому изоморфна сплетенному
произведению $S_{2}\wr S_{n}$. Множество изоморфизмов границы
$n$-мерного параллелепипеда по свойству сохранения или изменения
ориентации границы разбивается на два класса отображений, а
именно, на класс гомотопных и класс негомотопных тождественному
отображений. Принадлежность к тому или иному классу определяется
детерминантом матрицы представления изоморфизма, т.е. равенство
детерминанта матрицы единице означает принадлежность
соответствующего изоморфизма к классу гомотопных тождественному
отображений, а равенство $-1$ означает принадлежность к классу
негомотопных тождественному отображений. Откуда понятно, что
существование (отсутствие) гомотопного тождественному
антиподального отображения границы $n$-мерного параллелепипеда
определяется детерминантом диагональной матрицы, все элементы
которой равны $-1$, т.е. это свойство связано с четностью числа
$n$. В свою очередь, для всякого четного $n$ всегда существует
изоморфизм границы $n$-мерного параллелепипеда, квадрат которого
равен антиподальному отображению. Действительно, для этого
достаточно взять такую мономиальную подстановку базиса
$\mathbb{R}^{n}$, которая в произвольном разбиении базиса
$\{(e_{j})\}_{n}$ на непересекающиеся пары $(e_{i},e_{k})$
осуществляет вращение $(e_{i},e_{k})\mapsto(e_{k},-e_{i})$ всех
пар разбиения. С другой стороны, матрицы представления
изоморфизмов границы $n$-мерного параллелепипеда, квадрат которых
равен антиподальному отображению, принадлежат $SO(n)$ и им
соответствуют регулярные векторные поля соответствующей
нечетномерной сферы. Следовательно ортогональность пары таких
полей приводит к тому, что квадрат композиции соответствующих им
изоморфизмов границы параллелепипеда равен антиподальному
отображению. Таким образом, максимальное число попарно
ортогональных векторных полей нечетномерной сферы $S^{n-1}$ в
точности равно максимальному числу соответствующих изоморфизмов
границы $n$-мерного параллелепипеда.

Итак, пространство операторов, формирующих регулярные векторные
поля сферы мы будем искать в некоторой подгруппе $M(n)$ группы
мономиальных подстановок, состоящей исключительно из симметричных
и кососимметричных матриц, т.е. матриц, квадрат которых равен $E$
или $-E$. Пусть группа $M(2^{m})$ порождается $2m$-элементным
множеством $\{A_{1},\ldots,A_{m},B_{1},\ldots, B_{m}\}$, где
$A_{i}$--- блочно-диагональная матрица, состоящая из расположенных
по диагонали квадратных матриц порядка $k=2^{i}$, имеющих
ненулевые элементы везде кроме побочной диагонали, на которой
расположены $k$ единиц, а $B_{i}$ --- диагональная матрица,
диагональ которой заполнена последовательностью чередующихся
$2^{i-1}$ единиц и $2^{i-1}$ единиц со знаком минус, причем первые
$2^{i-1}$ единиц диагонали имеют знак плюс. Иначе говоря, оператор
$A_{i}$ действует на $\mathbb{R}^{n}$ так, что в каждом
подпространстве $\mathbb{R}^{k}$, на которые оно раскладывается в
прямую сумму, имеет место зеркальное отражение, т.е. перестановка
всех первых координат с последними, вторых с предпоследними, и так
далее. В свою очередь, оператор $B_{i}$ действует на
$\mathbb{R}^{n}$ так, что в каждом подпространстве
$\mathbb{R}^{k}$ имеет место зеркальная инверсия, т.е. изменение
знака второй половины координат. Тем самым, группу $M(2^{m})$
можно назвать {\it группой зеркальных симметрий}.

Как легко проверить, все генераторы группы $M(2^{m})$ являются
ортогональными симметричными операторами ($A_{i}^{2}=B_{i}^{2}=E$
для всех $i$), которые либо коммутируют либо антикоммутируют, а
именно, $A_{i}A_{j}=A_{j}A_{i}$, $B_{i}B_{j}=B_{j}B_{i}$ для всех
$i,j$, и $A_{i}B_{j}=-B_{j}A_{i}$ если $i\geq j$, но
$A_{i}B_{j}=B_{j}A_{i}$ если $i<j$. Пусть произвольный оператор
$C$, взятый из группы зеркальных симметрий, разложен в некоторое
произведение генераторов. Тогда, для того чтобы выяснить, является
ли $C$ симметричным или кососимметричным ортогональным оператором
нам необходимо вчислить чему равен квадрат оператора $C$, т.е.
установить равенство $C^{2}=E$ или $C^{2}=-E$. Прежде всего
заметим, что при перестановке произвольной пары генераторов
разложения оператора $C$ мы получим новый оператор $C'$, квадрат
которого по-прежнему равен квадрату оператора $C$, т.е.
$(C')^{2}=C^{2}$. Пусть в результате композиции перестановок
генераторов разложения оператора $C$ он приведен к виду $C'=AB$,
где ни в произведении $A$ ни в произведении $B$ нет генераторов с
повторяющимися индексами. Квадрат оператора $C'$ (равно как и $C$)
зависит от перестановочных свойств оператора $AB$. Действительно,
если $BA=AB$, то $(C')^{2}=E$, а если $BA=-AB$, то $(C')^{2}=-E$.
В свою очередь, перестановочные свойства оператора $AB$
определяются четностью числа $\delta$, равного количеству всех
возможных пар $A_{i}B_{j}$ в произведении $AB$, индексы которых
удовлетворяют условию $i\geq j$. Таким образом, в силу равенства
$BA=(-1)^{\delta}AB$ справедлива формула
\begin{equation}\label{delta}
  C^{2}=(AB)^{2}=(-1)^{\delta}E.
\end{equation}

Абстрактное представление группы зеркальных симметрий можно дать с
помощью $2m+1$ образующих и соотношений: $$M(2^{m})\simeq\langle
(a_{1},\ldots,a_{m}),(b_{1},\ldots,b_{m}),
  \varepsilon|a_{i}^{2}=b_{i}^{2}=\varepsilon^{2}=e,$$ $$a_{i}a_{j}=
  a_{j}a_{i},b_{i}b_{j}=b_{j}b_{i}, \varepsilon
  a_{i}=a_{i}\varepsilon,\varepsilon b_{i}=b_{i}\varepsilon,$$ $$
  a_{i}b_{j}=b_{j}a_{i} \forall i<j, a_{i}b_{j}=
  \varepsilon b_{j}a_{i} \forall i\geq j
  \rangle$$
Центр группы $M(2^{m})$ состоит из двух элементов
$\{e,\varepsilon\}$. Однако абстрактное задание группы зеркальных
симметрий не вполне удобно, и поэтому мы воспользуемся ее
представлением в виде неабелевой аддитивной группы. Пусть на
группе $\mathbb{Z}_{2}^{m}\times \mathbb{Z}_{2}^{m}$ задана
функция $\alpha$, принимающая значение 0 или 1 в соответствии с
отображением
\begin{equation}\label{alfa}
  \alpha :\mathbb{Z}_{2}^{m}\times \mathbb{Z}_{2}^{m}\rightarrow
  \mathbb{Z}_{2}:\alpha(a,b)=a_{1}b_{1}+a_{2}\sum_{1}^{2}b_{j}+\cdots
  +a_{i}\sum_{1}^{i}b_{j}+\cdots+a_{m}\sum_{1}^{m}b_{j},
\end{equation}
где $a_{i}$ это соответствующая компонента элемента $a$ а $b_{j}$
это соответствующая компонента $b$, и где подразумеается
суммирование единиц и нулей по модулю 2. Например, если $m=5$, то
$\alpha(a,b)=a_{1}b_{1}+ a_{2}(b_{1}+b_{2})+
a_{3}(b_{1}+b_{2}+b_{3})+ a_{4}(b_{1}+b_{2}+b_{3}+b_{4})+
a_{5}(b_{1}+b_{2}+ b_{3}+b_{4}+b_{5})$. Поскольку в сумму
$\alpha(a,b)$ включаются только такие равные единице произведения
$a_{i}b_{j}$, для которых выполнено условие $i\geq j$, то
$\alpha(a,b)=|\delta|\mod 2$, где число $\delta$ равно количеству
всех тех пар единиц элемента группы, которые принадлежат разным
компонентам произведения $\mathbb{Z}_{2}^{m}\times
\mathbb{Z}_{2}^{m}$ и чьи индексы удовлетворяют соотношению $i\geq
j$. Тем самым, если $\alpha(a,b)=0$, то в абстрактном предсталении
группы зеркальных симметрий соответствующее произведение
коммутирует, т.е. $ba=ab$, если же $\alpha(a,b)=1$, то
$ba=\varepsilon ab$. В свою очередь, если $\alpha(a',b)=0$, то
$aba'b'=aa'bb'$, если же $\alpha(a,b')=1$, то $aba'b'=\varepsilon
aa'bb'$. Следовательно группа зеркальных симметрий изоморфна
группе
$$\mathbb{Z}_{2}\leftthreetimes \left(\mathbb{Z}_{2}^{m}\times
\mathbb{Z}_{2}^{m}\right),$$ в которой произведение (сложение)
элементов $(s,a,b)$ и $(s',a',b')$ задано формулой
\begin{equation}\label{noabel}
  (s,a,b)+(s',a',b')=(s+s'+\alpha(a',b),a+a',b+b'),
\end{equation}
где отображение $\alpha$, конечно же, не является гомоморфизмом.
Далее, опираясь на абстрактное представление группы $M(2^{m})$
неабелевой аддитивной группой, мы будем всякий факторизованный по
центру ее элемент $ab$, соответствующий элементу
$(a,b)\in\mathbb{Z}_{2}^{m} \times \mathbb{Z}_{2}^{m}$, обозначать
символом $(I|J)$, где мультииндексы $I$, $J$ суть наборы индексов
единиц $a$ и $b$ соответственно. Таким образом, изучение группы
внутренних автоморфизмов группы зеркальных симметрий мы сводим к
изучению функции $\alpha$, определенной на группе
$\mathbb{Z}_{2}^{m}\times \mathbb{Z}_{2}^{m}$.

\section{Схема решения задачи}

Прежде всего установим основные свойства $\alpha$-функции.
\begin{enumerate}
  \item Если в формуле \ref{alfa} сначала перемножить все члены а
  затем сгруппировать слагаемые уже при $b_{i}$, то мы получим формулу
\begin{equation}\label{alfa2}
  \alpha(a,b)=b_{m}a_{m}+b_{m-1}\sum_{m-1}^{m}a_{j}+\cdots
  +b_{i}\sum_{m-i}^{m}a_{j}+\cdots+b_{1}\sum_{1}^{m}a_{j}
\end{equation}
  \item $\alpha(a,0)=\alpha(0,b)=0$
  \item $\alpha(a,b+a',b')=\alpha(a+a',b+b')=
  \alpha(a,b)+\alpha(a',b')+\alpha(a,b')+\alpha(a',b)$
  \item Поскольку при любом фиксированном ненулевом значении $a$
  функция $\alpha(a,b)$ равна сумме по модулю 2 некоторых (вполне
  определенных) компонентов $b$, то $\alpha(a,b)=1$ для всех тех
  $2^{m-1}$ значений $b$, для которых указанная сумма равна
  единице, и $\alpha(a,b)=0$ для всех тех $2^{m-1}$ значений $b$,
  для которых она равна нулю. Аналогично, при любом фиксированном
  ненулевом значении $b$ мы получим $2^{m-1}$ значений $a$, для
  которых $\alpha(a,b)=1$, и $2^{m-1}$ значений $a$, для
  которых $\alpha(a,b)=0$
  \item Так как мы имеем $2^{m}-1$ ненулевых значений $a$($b$), то
  функция $\alpha(a,b)$ принимает единичное значение на $2^{m-1}(2^{m}-1)$
  различных аргументах
  \item Если $\alpha(a,b)=\alpha(a',b')=1$, то равенства
  $\alpha(a,b+a',b')=1$ и $\alpha(a,b')+\alpha(a',b)=1$
  равносильны
  \item Если $(a,b)=(i|i)$, т.е. $a$ и $b$ равны элементам с
  единицей в $i$-ой компоненте и нулями во всех остальных
  компонентах, тогда $\alpha(i|i)=a_{i}b_{i}=1$
  \item Если $(i,j|i,j)=(i|i)+(j|j)$, где $i<j$, то $\alpha(i,j|i,j)=
  a_{i}b_{i}+a_{j}(b_{i}+b_{j})=1$
\end{enumerate}
Теперь наша задача сводится к тому, чтобы выбрать максимально
возможное множество $\{(a,b)\}_{p(m)}$ элементов
$\mathbb{Z}_{2}^{m}\times \mathbb{Z}_{2}^{m}$, для которых
выполняются условия $\alpha(a',b')=\alpha(a'',b'')=1$ и
$\alpha(a',b'+a'',b'')=1$, где $(a',b'),(a'',b'')\in
\{(a,b)\}_{p(m)}$. Однако, прежде чем обратиться к произвольному
$m$, воспользуемся основными свойствами $\alpha$-функции и решим
задачу в случае $m=1,2,3,4$.

Пусть $m=1$. Тогда в группе $\mathbb{Z}_{2}\times\mathbb{Z}_{2}$
имеется только один элемент $(1|1)$, на котором функция $\alpha$
принимает единичное значение. Следовательно для всякого вектора
$X=(X_{1},X_{2})$, лежащего на $S^{1}$, мы получим одно
ортогональное ему векторное поле $A_{1}B_{1}X=(-X_{2},X_{1})$.

Пусть $m=2$, где функция $\alpha$ принимает единичное значение на
шести своих аргументах. Представим функцию $\alpha$, определенную
на группе $\mathbb{Z}_{2}^{2} \times \mathbb{Z}_{2}^{2}$, таблицей
ее значений, лежащих на пересечении строк и столбцов, где столбцы
соответствуют элементам $b$ а строки -- элементам $a$, записанным
в виде двоичной записи номеров столбцов и строк (от 0 до 4),
причем старший разряд соответствует второму индексу а младший --
первому. Тогда легко получить три элемента, равные
$\{(1|1),(2|2),(1,2|1,2)\}$ (в двоичной форме они имеют вид
$\{(01,01),(10,10),(11,11)\}$), которые обдадают тем свойством,
что их попарные произведения (суммы) также отображаются функцией
$\alpha$ в единицу. Действительно, первые два элемента можно
выбрать в силу свойства 7 $\alpha$-функции, а третий элемент, c
учетом равенства $a+a=b+b=0$ и в силу свойства 2, не должен лежать
на строках и столбцах с номерами 0,1,2, но может лежать на
пересечении строки и столбца с номером 3. Но этот выбор следует
еще подтвердить. Для этого необходимо, последовательно складывая
элемент $(11,11)$ с уже выбранными элементами, проверить суммы
этих пар элементов на предмет равенства функции $\alpha$ от них
единице. В силу свойства 6, совместность третьего элемента с уже
выбранными элементами можно также проверить при помощи теста на
значение $\alpha$-функции в вершинах прямоугольника с двумя
диагонально расположенными вершинами в проверяемой паре элементов.
Таким образом, для всякого вектора $X=(X_{1},X_{2},X_{3},X_{4})$,
лежащего на $S^{3}$, мы получили три ортогональных ему векторных
поля $A_{1}B_{1}X=(-X_{2},X_{1},-X_{4},X_{3})$,
$A_{2}B_{2}X=(-X_{4},-X_{3},X_{2},X_{1})$,
$A_{1}A_{2}B_{1}B_{2}X=(-X_{3},X_{4},X_{1},-X_{2})$, которые
ортогональны друг другу.
\par\vspace{0.5cm}
{\begin{tabular}{|c|c|c|c|c|} \hline
     & 00 & 01 & 10 & 11 \\ \hline
  00 & 0 & 0 & 0 & 0 \\ \hline
  01 & 0 & {\bf 1} & 0 & 1 \\ \hline
  10 & 0 & 1 & {\bf 1} & 0 \\ \hline
  11 & 0 & 0 & 1 & {\bf 1} \\ \hline
\end{tabular}}
\par\vspace{0.5cm}

Пусть $m=3$, где функция $\alpha$, определенная на группе
$\mathbb{Z}_{2}^{3} \times \mathbb{Z}_{2}^{3}$, принимает
единичное значение на 28 своих аргументах. Зададим ее табличное
представление и в соответствии с ранее изложенной методикой найдем
элементы, совместные с задающими элементами
$\{(1|1),(2|2),(3|3)\}$ (в двоичной форме они имеют вид
$\{(001,001),(010,010),(100,100)\}$). Тем самым, мы получим
максимально возможное множество элементов, попарные произведения
(суммы) которых отображаются функцией $\alpha$ в 1, а именно, семь
элементов $\{(1|1),(2|2),(3|3),(1,2|1,2,3),(1,3|1,2),(2,3|1,3),
(1,2,3|2,3)\}$, в двоичной форме имеющих вид
$\{(001,001),(010,010),(100,100),$
$(011,111),(101,011),(110,101),$ $(111,110)\}$. Следовательно для
всякого вектора,
$X=(X_{1},X_{2},X_{3},X_{4},X_{5},X_{6},X_{7},X_{8})$, лежащего на
$S^{7}$, мы получим семь ортогональных ему векторных полей
$$A_{1}B_{1}X=(-X_{2},X_{1},-X_{4},X_{3},-X_{6},X_{5},-X_{8},X_{7}),$$
$$A_{2}B_{2}X=(-X_{4},-X_{3},X_{2},X_{1},-X_{8},-X_{7},X_{6},X_{5}),$$
$$A_{3}B_{3}X=(-X_{8},-X_{7},-X_{6},-X_{5},X_{4},X_{3},X_{2},X_{1}),$$
$$A_{1}A_{2}B_{1}B_{2}B_{3}X=(-X_{3},X_{4},X_{1},-X_{2},X_{7},-X_{8},-X_{5},X_{6})$$
$$A_{1}A_{3}B_{1}B_{2}X=(-X_{7},X_{5},X_{5},-X_{6},-X_{3},X_{4},X_{1},-X_{2}),$$
$$A_{2}A_{3}B_{1}B_{3}X=(-X_{5},X_{6},-X_{7},X_{8},X_{1},-X_{2},X_{3},-X_{4}),$$
$$A_{1}A_{2}A_{3}B_{2}B_{3}X=(-X_{6},-X_{5},X_{8},X_{7},X_{2},X_{1},-X_{4},-X_{3}),$$
которые ортогональны друг другу.
\par\vspace{0.5cm}
{\begin{tabular}{|c|c|c|c|c|c|c|c|c|} \hline
     & 000 & 001 & 010 & 011 & 100 & 101 & 110 & 111 \\ \hline
  000 & 0 & 0 & 0 & 0 & 0 & 0 & 0 & 0 \\ \hline
  001 & 0 & {\bf 1} & 0 & 1 & 0 & 1 & 0 & 1 \\ \hline
  010 & 0 & 1 & {\bf 1} & 0 & 0 & 1 & 1 & 0 \\ \hline
  011 & 0 & 0 & 1 & 1 & 0 & 0 & 1 & {\bf 1} \\ \hline
  100 & 0 & 1 & 1 & 0 & {\bf 1} & 0 & 0 & 1 \\ \hline
  101 & 0 & 0 & 1 & {\bf 1} & 1 & 1 & 0 & 0 \\ \hline
  110 & 0 & 0 & 0 & 0 & 1 & {\bf 1} & 1 & 1 \\ \hline
  111 & 0 & 1 & 0 & 1 & 1 & 0 & {\bf 1} & 0 \\ \hline
\end{tabular}}
\par\vspace{0.5cm}

Пусть $m=4$, где функция $\alpha$ принимает единичное значение на
120 элементах группы $\mathbb{Z}_{2}^{4}\times\mathbb{Z}_{2}^{4}$.
Если задать ее табличное представление, то в результате
последовательного поиска элементов, совместных с элементами
$\{(1|1),(2|2),(3|3),(4|4)\}$ (в двоичной форме они имеют вид
$\{(0001,0001),(0010,0010),(0100,0100),(1000,1000)\}$), мы получим
8 элементов, а именно, $\{(1|1),(2|2),(3|3),(4|4),(1,2|1,2,3),
(1,3|1,2,4),$ $(2,3|1,3,4),(1,2,3|2,3,4)\}$ (в двоичной форме они
имеют вид $\{(0001,0001),$ $(0010,0010),$ $(0100,0100),
(1000,1000),(0011,0111),(0101,1011),$ $(0110,1101),$
$(0111,1110)\}$), которые обладают тем свойством, что их попарные
произведения (суммы) также отображаются функцией $\alpha$ в
единицу. Кроме того, они обладают таким замечательным свойством,
что сумма по модулю 2 всех компонентов $b$ произвольного элемента
$(a,b)$ из множества $\{(a,b)\}_{p(4)}$ всегда равна единице. Имея
в виду соответствие элементов $(a,b)$ операторам $AB$, мы легко
получим восемь ортогональных друг другу векторных полей сферы
$S^{15}$.

Пусть теперь $m$ произвольно, причем $m=4d+c$, где $c=0,1,2,3$.
Случаи $m\leq 4$, в которых $d=0$ а $c=1,2,3$ и $d=1$ а $c=0$ мы
уже рассмотрели, поэтому пусть $m>4$. Тогда множество
$\{(a,b)\}_{p(m)}$ состоит из подмножества $\{(a,b)\}_{p(4)}$ и из
записанных в двоичной форме элементов $(a',b')$ вида
$$(2^{4k+1}a+2^{4k}\chi(a),2^{4k+1}b),$$ где $k=1,\ldots,d$;
дискретная функция $\chi(a)$ принимает нулевое значение в том
случае, когда сумма по модулю 2 всех компонентов $a$ равна
единице, и наоборот, если эта сумма равна нулю, то значение
функции $\chi(a)$ равно единице; $(a,b)\in\{(a,b)\}_{p(4)}$ для
$k=1,\ldots,d-1$ и $(a,b)\in\{(a,b)\}_{p(c)}$ для $k=d$, причем,
$p(0)=0$. Как легко проверить, $\alpha(a',b')=1$ а в силу того,
что сумма по модулю 2 всех компонентов $b$ и $a'$ по отдельности
равна единице, функция $\alpha$ от суммы двух произвольных
элементов множества $\{(a,b)\}_{p(m)}$ также равна единице, т.е.
множество $\{(a,b)\}_{p(m)}$ удовлетворяет всем необходимым
требованиям.

Итак, теперь мы умеем формировать на сфере $S^{n-1}$ систему из
$p(m)$ ортогональных векторных полей. Следовательно в каждой точке
сферы можно задать репер из $p(m)$ ортонормированных векторов, на
которое натягивается пространство $\mathbb{R}^{p(m)}$,
ортогональное вектору $X$. Всякое линейное преобразование
$\mathbb{R}^{p(m)}$, заданное обратимой (невырожденной) квадратной
матрицей порядка $p(m)$, трансформирует наш ортонормированный
репер произвольным образом, и тем самым переводит систему
ортогональных векторных полей в произвольную систему регулярных
линейно независимых векторных полей сферы $S^{n-1}$. Более того,
если элементы матрицы преобразования $\mathbb{R}^{p(m)}$ являются
произвольными непрерывными вещественнозначными функциями,
определенными на сфере, то мы получим произвольную систему
непрерывных линейно независимых векторных полей сферы $S^{n-1}$.
Заметим также, что перенести результат, полученный для $n=2^{m}$,
на случай $n=(2a+1)2^{m}$, не составляет труда, достаточно
разложить пространство $\mathbb{R}^{n}$ в прямую сумму $(2a+1)$
подпространств размерности $2^{m}$.

\appendix
\section{Псевдосфера}

Прежде всего заметим, что регулярные ортогональные векторные поля
3-мерной сферы индуцируют 3-мерное евклидово пространство.
Действительно, если на операторы $G_{1}=A_{1}B_{1}$,
$G_{2}=A_{2}B_{2}$, $G_{3}=A_{1}A_{2}B_{1}B_{2}$ натянуть
пространство $\langle G_{1},G_{2},G_{3}\rangle_{\mathbb{R}}$,
состоящее из элементов $X=x_{1}G_{1}+x_{2}G_{2}+x_{3}G_{3}$, то
$X^{2}=-(x_{1}^{2}+x_{2}^{2}+x_{3}^{2})E$. Следовательно тройка
генераторов $(G_{1},G_{2},G_{3})$ порождает классическую алгебру
Клиффорда, в которую вкладывается 3-мерное евклидово пространство
$\mathbb{R}^{3}$, состоящее из элементов $x=x_{1}e_{1}+
x_{2}e_{2}+x_{3}e_{3}$. Аналогично, регулярные ортогональные
векторные поля нечетномерной (n-1)-мерной сферы индуцируют
p(m)-мерное евклидово пространство. Однако обратимся к
псевдосферам, которые являются поверхностями уровня функции
$(x_{1})^{2}+\ldots + (x_{4})^{2}-(x_{5})^{2}-\ldots -
(x_{8})^{2}$.

Пусть в пространстве $\mathbb{R}^{8}$ задана метрика с сигнатурой
$(+,+,+,+,$ $-,-,-,-)$. Тогда касательное векторное поле
псевдосферы будет ортогональным (в этой метрике) к произвольному
вектору $X$, лежащему на ней. Для того чтобы сформировать
максимальное число ортогональных $X$ и ортогональных друг другу
регулярных векторных полей, обратимся к группе зеркальных
симметрий в случае $m=3$, откуда мы и возьмем четыре (максимально
возможное число) регулярных векторных поля, равных
$$G_{0}X=A_{3}X=(X_{8},X_{7},X_{6},X_{5},X_{4},X_{3},X_{2},X_{1}),$$
$$G_{1}X=A_{1}B_{1}X=(-X_{2},X_{1},-X_{4},X_{3},-X_{6},X_{5},-X_{8},X_{7}),$$
$$G_{2}X=A_{2}B_{2}X=(-X_{4},-X_{3},X_{2},X_{1},-X_{8},-X_{7},X_{6},X_{5}),$$
$$G_{3}X=A_{1}A_{2}B_{1}B_{2}X=(-X_{3},X_{4},X_{1},-X_{2},-X_{7},X_{8},X_{5},-X_{6}).$$
Легко заметить, что операторы $(G_{0},G_{1},G_{2},G_{3})$
составляют набор генераторов алгебры Клиффорда. Действительно,
если на них натянуть линейное пространство $\langle
G_{0},G_{1},G_{2},G_{3}\rangle_{\mathbb{R}}$, состоящее из
векторов $X=x_{0}G_{0}+x_{1}G_{1}+x_{2}G_{2}+x_{3}G_{3}$, то
$X^{2}=(x_{0}^{2}-x_{1}^{2}-x_{2}^{2}-x_{3}^{2})E$, а
следовательно мы получим 16-мерную алгебру, в которую вкладывается
псевдоевклидово пространство $\mathbb{R}^{4}$ с сигнатурой
$(+,-,-,-)$. Тем самым, мы фактически построили спинорное
представление псевдоортогональной группы $O(1,3)$. В самом деле,
пусть $g$ -- произвольный элемент группы внутренних автоморфизмов
пространства $\langle G_{0},G_{1},G_{2},G_{3}
\rangle_{\mathbb{R}}$, т.е. $gXg^{-1}=Y$, где
$Y=y_{0}G_{0}+y_{1}G_{1}+y_{2}G_{2}+y_{3}G_{3}$. Тогда
$Y^{2}=X^{2}$, откуда $y^{2}=x^{2}$, где $x,y\in\mathbb{R}^{4}$, и
поэтому всякому элементу $g$, действующему сопряжением на $\langle
G_{0},G_{1},G_{2},G_{3} \rangle_{\mathbb{R}}$, соответствует
некоторый псевдоортогональный оператор $A(g)$, действующий на
$\mathbb{R}^{4}$ так, что $A(g)x=y$, т.е. мы получили гомоморфизм
спинорной группы в псевдоортогональную группу, где отображение,
обратное этому гомоморфизму, и есть спинорное представление группы
$O(1,3)$.

Таким образом, регулярные векторные поля нашей псевдосферы вполне
могут служить источником пространственных симметрий макро- и
микрокосмоса.


\begin{thebibliography}{9}

\bibitem{LEN}С. Ленг, Алгебра, "Мир", Москва, 1968
\bibitem{VIN}Э. Б. Винберг, Курс алгебры, "Факториал", Москва,
1999
\bibitem{KUR}А. Г. Курош, Теория групп, Москва, 1967.
\bibitem{ADA}Adams, J. F. Vector fields on spheres. Ann. of Math. (2) 75 1962
603--632.
\bibitem{AMI}Ю.А. Аминов, Геометрия векторного поля, "Наука", Москва, 1990.

\end{thebibliography}
\end{document}